\documentclass[12pt]{amsart}
 \usepackage{amscd}
                     \usepackage{amssymb}
 \newcommand{\be}{\begin{equation}}
       \newcommand{\ee}{\end{equation}}
       \newcommand{\ba}{\begin{eqnarray}}
        \newcommand{\ea}{\end{eqnarray}}
 \newcommand{\ban}{\begin{eqnarray*}}
 \newcommand{\ean}{\end{eqnarray*}}

\def\XXint#1#2#3{{\setbox0=\hbox{$#1{#2#3}{\int}$}
     \vcenter{\hbox{$#2#3$}}\kern-.5\wd0}}

  \newcommand{\Pf}{\noindent {\bf Proof:} }
  \newcommand{\Rk}{\noindent {\bf Remark} }

 \newtheorem{theo}{Theorem}[section]

\begin{document}
 \newtheorem{defn}[theo]{Definition}
 \newtheorem{ques}[theo]{Question}
 \newtheorem{lem}[theo]{Lemma}
 \newtheorem{prop}[theo]{Proposition}
 \newtheorem{coro}[theo]{Corollary}
 \newtheorem{ex}[theo]{Example}
 \newtheorem{note}[theo]{Note}
 \newtheorem{conj}[theo]{Conjecture}
 \title[Shrinking gradient Ricci soliton]{Some rigidity results on shrinking gradient Ricci soliton}
 \author{Jianyu Ou}
 \email{oujianyu@xmu.edu.cn}
 \address{Department of Mathematics, Xiamen University}
 \author{Yuanyuan Qu}
\email{52285500012@stu.ecnu.edu.cn}
\address{School of Mathematical Sciences,Shanghai Key Laboratory of PMMP, East China Normal University, Shanghai 200241,
 China}
\author{Guoqiang Wu}
\email{gqwu@zstu.edu.cn}
\address{School of Science, Zhejiang Sci-Tech University, Hangzhou 310018, China}


 \date{}
 \maketitle

\begin{abstract} Suppose $(M^n, g, f)$ is a complete shrinking gradient Ricci soliton.
We give several rigidity results under some natural conditions, generalizing the results in  \cite{Petersen-Wylie,Guan-Lu-Xu}.  Using maximum principle, we prove that shrinking gradient Ricci soliton with constant scalar curvature $R=1$ is isometric to a finite quotient of $\mathbb{R}^2\times \mathbb{S}^2$, giving a new proof of the main results of Cheng-Zhou \cite{Cheng-Zhou}.
\end{abstract}

\section{Introduction}

For an $n$-dimensional complete Riemannian manifold $(M,g)$ and a smooth
potential function $f$ on $(M,g)$, the triple $(M, g, f)$ is called a
\emph{gradient shrinking Ricci soliton} or \emph{shrinker} 
if
\begin{align}\label{Eq1}
Ric+\mathrm{Hess}\,f=\frac 12g,
\end{align}
where $Ric$ is the Ricci curvature of $(M,g)$ and $\text{Hess}\,f$ is
the Hessian of $f$.  Shrinkers
are viewed as a natural extension of Einstein manifolds. More importantly,
shrinkers play an important role in the Ricci flow as they correspond to
some self-similar solutions and arise as limits of dilations of Type I
singularities in the Ricci flow. Shrinkers can also be
regarded as critical points of the Perelman's entropy functional and play
a significant role in Perelman's resolution of the Poincar\'e conjecture
\cite{Perelman1,Perelman2,Perelman3}.

 The study
of solitons has become increasingly important in both the study of the Ricci flow and metric measure space. Solitons play a direct
role as singularity dilations in the Ricci flow proof of uniformization.  In \cite{Perelman1}, Perelman introduced the ancient $\kappa$-solutions, which
play an important role in the singularity analysis, and he also proved that suitable blow down limit of ancient $\kappa$-solutions must be a shrinking gradient
Ricci soliton. In \cite{Perelman2}, Perelman proved that any two dimensional non-flat ancient $\kappa$-soluition must be the standard $S^2$, and he also classified three dimensional shrinking gradient Ricci soliton  under the assumption of nonnegative curvature and $\kappa$-noncollapseness. Due to the work of Perelman \cite{Perelman2}, Ni-Wallach \cite{Ni-Wallach}, Cao-Chen-Zhu \cite{Cao-Chen-Zhu}, the classification of three dimensional shrinking gradient Ricci soliton is complete. For more work on the classification of gradient Ricci soliton under various curvature condition, see \cite{Brendle1}, \cite{Cao-Chen}, \cite{Cao-Chen-Zhu}, \cite{Cao-Chen2}, \cite{Cao-Wang-Zhang}, \cite{Chen-Wang}, \cite{Eminenti-LaNave-Mantegazza}, \cite{Naber}, \cite{Petersen-Wylie}, \cite{Pigola-Rimoldi-Setti}, \cite{Wu-Wu-Wylie}, \cite{Zhang}.

In this paper, we study some rigidity problem about the shrinking gradient Ricci soliton.

In section 2, we provide some preliminary knowledge which will be used throughout the paper.

In section 3, we impose the additional assumption which is called condition A as follows:
\ban
|R(u, v, u, v)|\leq A \cdot Ric(u, u)
\ean
for any $|u|=|v|=1$ and $u\bot v$, where $A$ is a positive constant. There are so many examples satisfying condition A but having mixed sectional curvature. An explicit example is the K\"ahler shrinking gradient  Ricci soliton on $CP^2\# (- CP^2)$ constructed by Cao and Koiso independently.

Next  we can state the following  splitting result.
\begin{theo}\label{structure theorem} Let $(M^n, g, f)$ be a shrinking gradient Ricci soliton satisfying condition A, then the universal cover of $M$ is isometric to $\mathbb{R}^k\times \mathbb{N}^{n-k}$, where $\mathbb{N}$ is an $n-1$ dimensional shrinking gradient Ricci soliton with positive Ricci curvature.
\end{theo}
\Rk Petersen-Wylie \cite{Petersen-Wylie} and Guan-Lu-Xu \cite{Guan-Lu-Xu} proved the above theorem independently if $(M^n, g, f)$ has nonnegative sectional curvature.

In section 4,
at first we define a symmetric two tensor $h$ by   $h(u, v)=\sum_{i, j=1}^n R(u, e_i, v, e_j)Ric(e_i, e_j)$, where $\{e_i\}_{i=1}^n$ are local orthonormal basis, then we state our main results as follows.
\begin{theo} Let $(M^n, g, f)$ be a shrinking gradient Ricci soliton, if $Ric\geq 0$ and $h\leq \frac{1}{2}Ric$, then the universal cover of $M$ is isometric to $\mathbb{R}^{k}\times \mathbb{N}^{n-k}$, where $\mathbb{N}$ is a compact Einstein manifold.
\end{theo}
Based on the above theorem, we give a new proof of Theorem 1.4 in \cite{Petersen-Wylie2}.
\begin{coro}\cite{Petersen-Wylie2} If $(M^n, g, f)$ is a shrinking gradient Ricci soliton with nonnegative sectional curvature and $R\leq 1$, then the universal cover of $M$ is isometric to either $\mathbb{R}^n$ or $\mathbb{S}^2\times \mathbb{R}^{n-2}$
\end{coro}

In section 5,  we focus our attention on $4$-dimensional gradient shrinking Ricci solitons with constant scalar curvature. Recall that in Petersen and Wylie's paper \cite{Petersen-Wylie}, a gradient Ricci soliton $(M, g)$ is said to be rigid if it is isometric to a quotient $N \times \mathbb{R}^k$, the product soliton of an Einstein manifold $N$ of positive scalar curvature with the Gaussian soliton $\mathbb{R}^k$.
Therefore, a gradient Ricci soliton has constant scalar curvature if it is rigid.
Conversely,  for the complete shrinking case, Prof. Huai-Dong Cao raised the following

\smallskip	
\noindent {\bf Conjecture}:
Let $(M^n, g, f)$, $n\geq 4$, be a complete $n$-dimensional gradient shrinking Ricci soliton. If $(M, g)$ has constant scalar curvature, then it must be rigid,
i.e., a finite quotient of $\mathbb{N}^k\times \mathbb{R}^{n-k}$ for some Einstein manifold $\mathbb{N}$ of positive scalar curvature.

\smallskip	

About this conjecture, Petersen-Wylie \cite{Petersen-Wylie2} proved the following interesting result.
\begin{theo}[\cite{Petersen-Wylie2}] A Ricci shrinker is rigid iff it has constant scalar curvature and is radially flat, that is, the sectional curvautre
\begin{align*}
K(\nabla f, \cdot)=0.
\end{align*}
\end{theo}
Later, Fern\'{a}ndez-L\'{o}pez and Garc\'\i a-R\'\i o \cite{FR16} characterize the rigidity using the rank of Ricci curvature.
\begin{theo}[\cite{FR16}]\label{constant rank} A  Ricci shrinker is rigid iff it has constant scalar curvature and the Ricci curvature has constant rank.
\end{theo}
In the same paper, they also proved that the possible value of $R$ is $\{0, 1, ..., \frac{n-1}{2}, \frac{n}{2}\}$.
In dimension $4$,   If $R=0$, $(M^4, g)$ is isometric to $\mathbb{R}^4$;  if $R=\frac{3}{2}$, then $(M^4, g)$ is  isometric to $\mathbb{R}\times \mathbb{S}^3$; if $R=2$, then $(M^4, g)$ is  isometric to a compact Einstein manifold with $Ric=\frac{1}{2}g$.

Recently, Cheng-Zhou \cite{Cheng-Zhou} proved a four dimensional Ricci shrinker with $R=1$ is rigid. They applied $\Delta_f$ to the  quantity
\begin{align*}
tr(Ric^3)-\frac{1}{4},
\end{align*}
and they got the following nice inequality
\begin{align*}
\Delta_f  \left(f( tr(Ric^3)-\frac{1}{4})\right)\geq 9 f( tr(Ric^3)-\frac{1}{4}),
\end{align*}
at last they used the integration by parts to derive that
\begin{align*}
tr(Ric^3)-\frac{1}{4}=0
\end{align*}
over $M$, and this implies that $\lambda_1+\lambda_2=0$, so the Ricci curvature has rank $2$, finally they obtain the rigidity by Theorem \ref{constant rank}.

We want to point remark that  $(3.11)$ in Cheng-Zhou \cite{Cheng-Zhou} gives that
\begin{align*}
\frac{1}{3}tr(Ric^3)-\frac{1}{12}=\lambda_2\lambda_3\lambda_4,
\end{align*}
this implies that the quantity they used is $\sigma_3(Ric)$, since $\lambda_1=0$ in this situation.

We restate the main theorem in \cite{Cheng-Zhou} as follows.
\begin{theo}[\cite{Cheng-Zhou}] Let $(M, g, f)$ be a 4-dimensional complete noncompact shrinking gradient Ricci soliton. If $M$ has constant scalar curvature $R=1$, then it must be isometric to a finite quotient of $\mathbb{R}^2\times S^2$.
\end{theo}
Our new proof is inspired by  \cite{Petersen-Wylie},  where they assumed the sectional curvature is nonnegative. Denote the eigenvalues of Ricci curvature by $\lambda_1\leq \lambda_2\leq \cdot\cdot\cdot\leq \lambda_n$.
  then it is easy to observe that $Rm*Ric\geq 0$. They applied $\Delta_f$ directly to the sum of the smallest $k$ eigenvalues and obtained
\begin{align*}
\Delta_f(\lambda_1+\lambda_2+ \cdot\cdot\cdot+ \lambda_k) \leq (\lambda_1+\lambda_2+\cdot\cdot\cdot+  \lambda_k)
\end{align*}
holds in the barrier sense, at last they derived that the Ricci curvature has constant rank by standard maximum principle.

In the setting of $R=1$, it is not that easy, and the most important thing is to derive the following inequality,
 \ban
 &&\Delta_f (\lambda_1+\lambda_2)\\
 &&\leq \frac{2\nabla f\cdot \nabla(\lambda_1+\lambda_2)+(\lambda_1+\lambda_2)-2(\lambda_1^2+\lambda_2^2)}{f}(1-2\lambda_1-2\lambda_2)\\
 &&\quad\quad -2(\lambda_1+\lambda_2)+4(\lambda_1^2+\lambda_2^2),
 \ean
then
\ban
 \frac{\lambda_1+\lambda_2}{f}
\ean
satisfies
\ban
\Delta_f\frac{\lambda_1+\lambda_2}{f}\leq -0.9 \cdot\frac{\lambda_1+\lambda_2}{f}
\ean
outside a compact set, next we can use similar trick as \cite{Munteanu-Wang4} to obtain a uniform positive lower bound of $\lambda_1+\lambda_2$, this is impossible unless $\lambda_1+\lambda_1\equiv 0$.


{\bf Acknowledgments}.
The authors would like to thank Professor Mijia Lai,  Professor Yu Li, Professor Xiaolong Li and Professor Fengjiang Li  for their helpful discussions. Wu is also grateful to Professor Xi-Nan Ma for his constant encouragement.

\section{Preliminary}

Suppose $(M^n, g, f)$ is a shrinking gradient Ricci soliton $\nabla^2 f+Ric=\frac{1}{2}g$.
At first we recall some basic formulas which will be used during the paper.

\ba
&&d R=2 Ric(\nabla f),\\ \label{second bianchi}
&& R+\Delta f=\frac{n}{2},\\
&& R+|\nabla f|^2=f,\\  \label{f he tidu}
&& \Delta_f R=R-2|Ric|^2,\\
&& \Delta_f R_{ij}=R_{ij}-2 R_{ikjl}R_{kl},\label{elliptic equation}
\ea
where  $\Delta_f Ric=\Delta Ric -\nabla_{\nabla f}Ric$  in the last formula.

Next we state the estimate of potential function $f$ in Cao-Zhou \cite{Cao-Zhou}.
\begin{theo}[\cite{Cao-Zhou}]\label{potential estimate}  Suppose $(M^n, g, f)$ is an noncompact shrinking gradient Ricci soliton, then there exist $C_1$ and $C_2$ such that
\ban
\left(\frac{1}{2}d(x, p)-C_1\right)^2\leq f(x)\leq \left(\frac{1}{2}d(x, p)+C_2\right)^2,
\ean
where $p$ is the minimal point of  $f$.
\end{theo}
\Rk. Chen \cite{Chen} proved that any shrinking gradient Ricci soliton has $R\geq 0$, so  due to $R+|\nabla f|^2=f$ we derive that $|\nabla f|(x)\leq \frac{1}{2}d(x, p)+C_2$.

The following splitting theorem for shrinking gradient Ricci soliton will be important for us.
\begin{theo}[Naber, \cite{Naber}]\label{naber's splitting theorem}For any $n$-dimensional shrinking gradient Ricci soliton $(M^n, g, f)$  with bounded curvature and a sequence of points $x_i\in M$ going to infinity along an integral curve of $\nabla f$, by choosing a subsequence if necessary, $(M, g, x_i)$ converges smoothly to a product manifold $\mathbb{R}\times \mathbb{N}^{n-1}$,  where $\mathbb{N}$ is a shrinking gradient Ricci soliton.
\end{theo}

%
\section{Structure of shrinking gradient Ricci soliton satisfying condition A}

There is one explicit curvature condition called 2-nonnegative flag curvature which is defined as follows.
\begin{defn}
$(M^n, g)$ is said to have 2-nonnegative flag curvature if
\ban
R(e_1, e_2, e_1, e_2)+R(e_1, e_3, e_1, e_3)\geq 0
\ean
for any othonormal basis $e_1, e_2, e_3$.
\end{defn}

\Rk Qu-Wu \cite{Qu-Wu} proved that 2-nonnegative flag curvature implies condition A. 2-nonnegative flag curvature is also considered in Li-Ni \cite{Li-Ni}.

 \begin{prop}\label{constant rank} Suppose $(M^n, g, f)$ is a shrinking gradient Ricci soliton satisfying condition A, then the rank of Ricci curvature tenor is constant.
 \end{prop}
 \Pf Denote the eigenvalues of Ricci curvature by $\lambda_1\leq \lambda_2\leq \cdots\leq \lambda_n$.

 \noindent\textbf{Claim.} $\Delta_f \lambda_1\leq (1+2 A (n-1)) \lambda_1$ in the barrier sense.

Actually, at $x$, assume $e_1$ satisfies  $Ric(x)(e_1, e_1)=\lambda_1(x)$, then extend $e_1$ to an orthonormal basis $\{e_1, e_2, \cdots, e_n\}$ such that $\{e_i\}_{i=1}^n$ are the eigenvectors of $Ric(x)$ corresponding to eigenvalues $\{\lambda_1, \lambda_2, \cdots, \lambda_n\}$.  Take parallel transport of $e_1$ along all the geodesics from $x$, then in a neighborhood $B(x, \delta)$ we get a smooth function $Ric(y)(e_1(y), e_1(y))$ satisfying $Ric(y)(e_1(y), e_1(y))\geq \lambda(y)$ and  $Ric(x)(e_1(x), e_1(x))= \lambda_1(x)$.

So at $x\in M$,
\ban
&&\Delta_f Ric(e_1, e_1)\\
&&=(\Delta_f Ric)(e_1, e_1)\\
&&=Ric(e_1, e_1)-2\sum_{i=1}^n R(e_1, e_i, e_1, e_i)Ric(e_i, e_i)\\
&&\leq Ric(e_1, e_1)+2 A (n-1)Ric(e_1, e_1)|Ric|\\
&&=Ric(e_1, e_1)\left(1+2 A (n-1)|Ric|  \right).
\ean\qed

Suppose there exists $q\in M$ such that $\lambda_1(q)=0$, then by the strong maximum principle, we get $\lambda_1\equiv 0$ on $M$.

Similar argument implies that $\lambda_1+\lambda_2+\cdots+\lambda_k$ satisfies
\ban
\Delta_f(\lambda_1+\lambda_2+\cdots+\lambda_k)\leq (1+2A (n-1)|Ric|)(\lambda_1+\lambda_2+\cdots+\lambda_k)
\ean
in the barrier sense for any $2\leq k\leq n$. So we can again use the maximum principle to derive that either $\lambda_1+\lambda_2+\cdots+\lambda_k>0$ or $\lambda_1+\lambda_2+\cdots+\lambda_k \equiv 0$ on $M$. \qed
\begin{prop}\label{parallel translation} Under the same assumption with Proposition \ref{constant rank}, the kernel of Ricci curvature is a parallel subbundle.
\end{prop}

\Pf Given any section $V\in ker(Ric)$, choose local orthonormal basis $\{e_1, e_2, \cdots, e_n\}$. Due to the nonnegativity of Ricci curvature, $Ric(V, V)=0$ is the same as $Ric(V)=0$. According to condition A, $Ric(V, V)=0$ also implies $R(V, e_i, V, e_i)=0$.
\ban
&&\Delta_f Ric(V, V)\\
&&=\Delta Ric(V, V)-\nabla_{\nabla f}Ric(V, V)\\
&&=\nabla_k\nabla_k(R_{ij}V^i V^j)-(\nabla_{\nabla f}Ric)(V, V)-2 Ric(\nabla_{\nabla f}V, V)\\
&&=(\Delta Ric)(V, V)+4\nabla_k R_{ij}\nabla_k V^i V^j+2Ric(\Delta V, V)+2 R_{ij}\nabla_k V^i \nabla_k V^j\\
&&\  -(\nabla_{\nabla f}Ric)(V, V)-2 Ric(\nabla_{\nabla f}V, V)\\
&&=(\Delta_f Ric)(V, V)+4\nabla_k R_{ij}\nabla_k V^i V^j+2Ric(\Delta_f V, V)+2 R_{ij}\nabla_k V^i \nabla_k V^j\\
&&=Ric(V, V)-2R(V, e_i, V, e_i)Ric(e_i, e_i)+4\nabla_k R_{ij}\nabla_k V^i V^j+2 R_{ij}\nabla_k V^i \nabla_k V^j\\
&&=4\nabla_k R_{ij}\nabla_k V^i V^j+2 R_{ij}\nabla_k V^i \nabla_k V^j.
\ean
Since
\ban
&&\nabla_k(R_{ij}\nabla_k V^i V^j)=\nabla_k R_{ij}\nabla_k V^i V^j +R_{ij}\Delta V^i V^j +R_{ij}\nabla_k V^i \nabla_k V^j\\
&&\quad =\nabla_k R_{ij}\nabla_k V^i V^j  +R_{ij}\nabla_k V^i \nabla_k V^j,
\ean
so
\ban
0=\Delta_f Ric(V, V)=-2 R_{ij}\nabla_k V^i \nabla_k V^j,
\ean
hence $\nabla_k V\in ker(Ric)$. \qed

Combining Proposition \ref{constant rank}, Proposition \ref{parallel translation} with De Rham's splitting Theorem, we can get the following structure result for shrinking gradient Ricci soliton.

\begin{theo}\label{structure theorem} Let $(M^n, g, f)$ be a shrinking gradient Ricci soliton satisfying condition A, then the universal cover of $M$ is isometric to $\mathbb{R}^k\times \mathbb{N}^{n-k}$, where $\mathbb{N}$ has positive Ricci curvature.
\end{theo}

Assume $(M^n, g, f)$ has nonnegative sectional curvature, Petersen-Wylie \cite{Petersen-Wylie} and Guan-Lu-Xu \cite{Guan-Lu-Xu} proved the above Theorem independently. Moreover, Munteanu-Wang \cite{Munteanu-Wang4} obtained that $(M^n, g, f)$ is compact if it has nonnegative sectional curvature and positive Ricci curvature. So $\mathbb{N}$ is compact in the above Theorem.
 Here our condition A is weaker than theirs. In \cite{Qu-Wu,Wu-Wu}, the authors derived  the soliton $(M^4, g, f)$ is also compact under condition A and other natural conditions.

\section{Rigidity via  nonnegative sectional curvature}

Define a symmetric two tensor $h$ by   $h(u, v)=\sum_{i, j=1}^n R(u, e_i, v, e_j)Ric(e_i, e_j)$, where $\{e_i\}_{i=1}^n$ are local orthonormal basis.

\begin{theo} Let $(M^n, g, f)$ be a shrinking gradient Ricci soliton, if $Ric\geq 0$ and $h\leq \frac{1}{2}Ric$, then the universal cover of $M$ is isometric to $\mathbb{R}^{k}\times \mathbb{N}^{n-k}$, where $\mathbb{N}$ is a compact Einstein manifold.
\end{theo}
\Pf Because $Ric\geq 0$ and $h\leq \frac{1}{2}Ric$,
\ban
\langle h-\frac{1}{2}Ric, Ric\rangle=R_{ikjl}R_{ij}R_{kl}-\frac{1}{2}|Ric|^2\leq 0.
\ean
Recall the formula (\ref{elliptic equation}), direct calculation gives
\ban
\frac{1}{2}\Delta_f |Ric|^2=|\nabla Ric|^2+ |Ric|^2-2 R_{ikjl}R_{ij}R_{kl}
\ean
\noindent\textbf{Claim.} $|\nabla Ric|\equiv 0$ on $M$.

Actually, integrating the above identity over $B(p, r)$, the right hand is always nonnegative, so the goal is to prove
\ban
\lim_{r\rightarrow\infty}\int_{B(p, r)}\Delta_f |Ric|^2 e^{-f}=0.
\ean
Thanks to the estimate $\int_{M}\left(|Ric|^2+|\nabla Ric|^2\right) e^{-f}<C$ in \cite{Li-Wang}, it is easy to see that there exists a sequence
$r_i\rightarrow \infty$ such that $\int_{\partial B(p, r_i)}\left(|Ric|^2+|\nabla Ric|^2\right) e^{-f}\rightarrow 0$. Hence
\ban
&&\left|\int_{B(p, r_i)}\Delta_f |Ric|^2 e^{-f}\right|\\
&&=\left|\int_{\partial B(p, r_i)}\langle\nabla |Ric|^2 e^{-f}, \nu\rangle\right|\\
&&\leq \int_{\partial B(p, r_i)}2 |Ric||\nabla Ric|\ e^{-f}\\
&&\leq \int_{\partial B(p, r_i)} \left(|Ric|^2+|\nabla Ric|^2\right) e^{-f}\rightarrow 0.\qed
\ean
So $|\nabla Ric|\equiv 0$ on $M$, De Rham's splitting Theorem implies that the universal cover of $M$ is isometric to $\mathbb{R}^{k}\times \mathbb{N}^{n-k}$, where $\mathbb{N}$ is a compact Einstein manifold. \qed

\begin{coro}\label{three proof} If $(M^n, g, f)$ is a shrinking gradient Ricci soliton with nonnegative sectional curvature and $R\leq 1$, then the universal cover of $M$ is isometric to either $\mathbb{R}^n$ or $\mathbb{S}^2\times \mathbb{R}^{n-2}$
\end{coro}
\Pf Choose local orthonormal basis $\{e_i\}_{i=1}^n$, due to the assumption, for any $i$,
\ban
2 Ric(e_i, e_i)=Ric(e_i, e_i)+\sum_{j\neq i}R(e_i, e_j, e_i, e_j)\leq Ric(e_i, e_i)+\sum_{j\neq i} Ric(e_j, e_j)=R\leq 1,
\ean
hence $Ric(e_i, e_i)\leq \frac{1}{2}$,    i.e. $Ric\leq \frac{1}{2}g$.
So
\ban
h(u, u)=\sum_{i=1}^n R(u, e_i, u, e_i)Ric(e_i, e_i)\leq \frac{1}{2}\sum_{i=1}^n R(u, e_i, u, e_i) =\frac{1}{2}Ric(u, u).
\ean
 Then from the above Theorem we know that the universal cover of $M$ is isometric to $\mathbb{R}^{k}\times \mathbb{N}^{n-k}$, where $\mathbb{N}$ is a compact Einstein manifold. Since $R\leq 1$, $\mathbb{N}$ has to be $\mathbb{S}^2$.
 \qed

\Rk. Corollary 4.2 appeared in  Petersen-Wylie's paper \cite{Petersen-Wylie2}, but our proof is different from theirs.  Under the same assumption, they first got the scalar curvature is identically one using Naber's result \cite{Naber}, then the conclusion follows from their main  theorem that  shrinking gradient soliton satisfying $0\leq Ric\leq \frac{1}{2} g$ and constant scalar curvature condition must be rigid.

According to the above discussion, we can give "the third  proof" using our structure Theorem \ref{structure theorem},

\noindent\textbf{Third proof of Corollary \ref{three proof}}

   Obviously the condition A holds, then the universal cover of $M^n$ is isometric to $\mathbb{R}^k\times \mathbb{N}^{n-k}$ by Theorem \ref{structure theorem}, where $\mathbb{N}$ is a shrinking gradient Ricci soliton with positive Ricci curvature. Because Naber's theorem implies the scalar curvature is identically one, $\nabla f$ is identically zero using formula (\ref{second bianchi}), i.e. $Ric=\frac{1}{2}g$ on $\mathbb{N}$, so $\mathbb{N}$ is $\mathbb{S}^2$.\qed
\begin{theo} Let $(M^n, g, f)$ be a shrinking gradient Ricci soliton with bounded curvature, if   $Ric\geq 0$ and  the scalar curvature $R<1-\delta$ for some $0<\delta<1$, then $M^n$ is flat.
\end{theo}
\Pf Suppose on the contrary, then the strong maximum principle gives  $R>0$ on $M$.

Because the curvature is bounded,  formula (\ref{f he tidu}) and quadratic growth of $f$ from Theorem \ref{potential estimate}, it is easy to see that outside a compact set, $f$ has no critical point. Formula (\ref{second bianchi}) implies that the scalar curvature is always increasing along the integral curve of $f$, hence the scalar curvature has a positive lower bound. Next we can apply Theorem \ref{naber's splitting theorem} to obtain that $(M, g, f)$ converge along the integral curve of $\nabla f$ to $\mathbb{R}\times \mathbb{N}^{n-1}$, where $\mathbb{N}^{n-1}$ is a nonflat shrinking gradient Ricci soliton satisfying the same assumption, then we play the same game on $\mathbb{N}$. When $\mathbb{N}$ is of dimension 4, According to the main Theorem in Munteanu-Wang \cite{Munteanu-Wang5} we get  the asymptotic limit of $\mathbb{N}$ is either  $\mathbb{R}\times \mathbb{S}^3$ or $\mathbb{R}^2\times \mathbb{S}^2$ or  their quotients. In any case the scalar curvature of the asymptotic limit can't be smaller than $1-\delta$. Contradiction.     \qed

\Rk In dimension 4, Munteanu-Wang \cite{Munteanu-Wang3} proved that bounded scalar curvature implies bounded curvature operator, so the bounded curvature assumption can be removed.

%
%
%
%

\section{constant scalar curvature}

To prove the main theorem, it is necessary to derive the nonnegativity of Ricci curvature. Actually,
\begin{lem}[\cite{FR16},\cite{Cheng-Zhou}] Suppose $(M^4, g, f)$ is a four dimensional Ricci shrinker with $R=1$, then $Ric\geq 0$.
\end{lem}

The following formula which is derived in plays an important role.
\ba\label{5.6}
\nabla_{\nabla f}Ric =Ric\circ(Ric-\frac{1}{2}g)+R(\nabla f, \cdot, \nabla f, \cdot).
\ea

Next we state the main result proved in Cheng-Zhou \cite{Cheng-Zhou} and give a new proof.
\begin{theo}[\cite{Cheng-Zhou}] Let $(M, g, f)$ be a 4-dimensional complete noncompact shrinking gradient Ricci soliton. If $M$ has constant scalar curvature $R=1$, then it must be isometric to a finite quotient of $\mathbb{R}^2\times S^2$.
\end{theo}

\Pf
Denote the eigenvalues of Ricci curvature by $\lambda_1\leq \lambda_2\leq  \lambda_3 \leq \lambda_4$.

 \noindent\textbf{Claim.}
 \ban
 &&\Delta_f (\lambda_1+\lambda_2)\\
 &&\leq \frac{2\nabla f\cdot \nabla(\lambda_1+\lambda_2)+(\lambda_1+\lambda_2)-2(\lambda_1^2+\lambda_2^2)}{f}(1-2\lambda_1-2\lambda_2)\\
 &&\quad\quad -2(\lambda_1+\lambda_2)+4(\lambda_1^2+\lambda_2^2)
 \ean
 in the barrier sense.

Actually, at $x$, because $R=1$, $Ric(\nabla f)=0$, so we choose $e_1=\frac{\nabla f}{|\nabla f|}$, then extend $e_1$ to an orthonormal basis $\{e_1, e_2, e_3, e_4\}$ such that $\{e_i\}_{i=1}^4$ are the eigenvectors of $Ric(x)$ corresponding to eigenvalues $\{\lambda_1, \lambda_2, \lambda_3, \lambda_4\}$.  Take parallel transport of $\{e_i\}_{i=1}^4$ along all the geodesics from $x$, then in a neighborhood $B(x, \delta)$ we get a smooth function $u(y)=Ric(y)(e_1(y), e_1(y))+Ric(y)(e_2(y), e_2(y))$ satisfying $u(y)\geq \lambda_1(y)+\lambda_2(y)$ and  $u(x)= \lambda_1(x)+\lambda_2(x)$.

So, at $x$, by (\ref{5.6}),
\ban
&&-2R(\nabla f, e_i, \nabla f, e_i)=2\nabla_k f_{ii}f_k +2 \lambda_i^2-\lambda_i\\
&&=-2(\nabla_{\nabla f} Ric)(e_i, e_i)+2 \lambda_i^2-\lambda_i,
\ean

Hence
\ban
2R(e_1, e_2, e_1, e_2)=\frac{2\nabla f\cdot \nabla u +u-2 u^2}{f}.
\ean

\ban
&&\Delta_f u (x)=\Delta_f\left( Ric(y)(e_1(y), e_1(y))+Ric(y)(e_2(y), e_2(y))\right)|_{y=x}\\
&&=(\Delta_f Ric )(e_1, e_1)(x)+(\Delta_f Ric )(e_2, e_2)(x)\\
&&= \lambda_1+\lambda_2-2\sum_{i=1}^4 R_{1i1i}\lambda_i\\
&&=(\lambda_1+\lambda_2)-2R_{1212}(\lambda_1+\lambda_2)-2\lambda_3(\lambda_3-R_{3434})-2\lambda_4(\lambda_4-R_{3434})\\
&&=(\lambda_1+\lambda_2)-2R_{1212}(\lambda_1+\lambda_2)-2(\lambda_3^2+\lambda_4^2)+2(\lambda_3+\lambda_4)R_{3434}\\
&&=(\lambda_1+\lambda_2)-2R_{1212}(\lambda_1+\lambda_2)-2(\lambda_3^2+\lambda_4^2)+\left(2R_{1212}+(\lambda_3+\lambda_4)-(\lambda_1+\lambda_2)\right)(\lambda_3+\lambda_4)\\
&&=(\lambda_1+\lambda_2)-(\lambda_1+\lambda_2)(\lambda_3+\lambda_4)+2R_{1212}(\lambda_3+\lambda_4-\lambda_1-\lambda_2)+(\lambda_3+\lambda_4)^2-2(\lambda_3^2+\lambda_4^2)\\
&&=2R_{1212}(1-2\lambda_1-2\lambda_2)+(\lambda_1+\lambda_2)+(\lambda_1+\lambda_2)(1-\lambda_3-\lambda_4)+(\lambda_3+\lambda_4)^2-2(\lambda_3^2+\lambda_4^2)\\
&&=\frac{2\nabla f\cdot \nabla(\lambda_1+\lambda_2)+(\lambda_1+\lambda_2)-2(\lambda_1^2+\lambda_2^2)}{f}(1-2\lambda_1-2\lambda_2)\\
&&+(\lambda_1+\lambda_2)^2+(1-\lambda_1-\lambda_2)^2-2(\frac{1}{2}-\lambda_1^2-\lambda_2^2)\\
&&=\frac{2\nabla f\cdot \nabla(\lambda_1+\lambda_2)+(\lambda_1+\lambda_2)-2(\lambda_1^2+\lambda_2^2)}{f}(1-2\lambda_1-2\lambda_2)\\
&&+(\lambda_1+\lambda_2)^2+1-2(\lambda_1+\lambda_2)+(\lambda_1+\lambda_2)^2-1+2(\lambda_1+\lambda_2)^2\\
&&=\frac{2\nabla f\cdot \nabla(\lambda_1+\lambda_2)+(\lambda_1+\lambda_2)-2(\lambda_1^2+\lambda_2^2)}{f}(1-2\lambda_1-2\lambda_2)-2(\lambda_1+\lambda_2)+4(\lambda_1^2+\lambda_2^2)\\
&&=\frac{2\nabla f\cdot \nabla u+u-2u^2}{f}(1-2u)-2 u+4u^2,
\ean
this finish the proof of the Claim.

Next we consider the function $\frac{\lambda_1+\lambda_2}{f}$,

\begin{align*}
&\Delta_f\frac{\lambda_1+\lambda_2}{f}=\frac{1}{f}\Delta_f(\lambda_1+\lambda_2)+\Delta_f{\frac{1}{f}}(\lambda_1+\lambda_2)+2\nabla \frac{1}{f}\cdot \nabla(\lambda_1+\lambda_2)\\
&\leq \frac{1}{f} \left\lbrace\frac{2\nabla f\cdot \nabla(\lambda_1+\lambda_2)+(\lambda_1+\lambda_2)-2(\lambda_1^2+\lambda_2^2)}{f}(1-2\lambda_1-2\lambda_2)\right.\\ &\quad\quad\quad\left.-2(\lambda_1+\lambda_2)+4(\lambda_1+\lambda_2)^2   \right\rbrace +\left( \frac{1}{f}+\frac{1}{f^2}   \right)(\lambda_1+\lambda_2)-\frac{2\nabla f\cdot\nabla(\lambda_1+\lambda_2)}{f^2}\\
&=\frac{\lambda_1+\lambda_2}{f}\left\lbrace-2+4(\lambda_1+\lambda_2)+\frac{1-2(\lambda_1+\lambda_2)}{f}  (1-2\lambda_1-2\lambda_2)
\right.\\ &\quad\quad\quad\quad\quad\quad-\left.\frac{4\nabla f\cdot\nabla(\lambda_1+\lambda_2)}{f}+1+\frac{1}{f}   \right\rbrace
\end{align*}

By Theorem \ref{naber's splitting theorem}, we know the asymptotic limit of $(M^4, g)$ along the integral curve of $f$ is $\mathbb{R}\times \mathbb{N}^3$, where $\mathbb{N}^3$ is a three dimensional shrinking gradient Ricci soliton, hence is a finite quotient of $\mathbb{R}\times \mathbb{S}^2$ due to $R=1$. So $\lambda_1+\lambda_2$ tends to zero at infinity.

Because $R=1$, Munteanu-Wang \cite{Munteanu-Wang3} obtained that Riemannian curvature is bounded, hence its derivative is also bounded due to Shi's estimate. So $|\nabla f \cdot \nabla(\lambda_1+\lambda_2)|$ is bouned by $C\cdot |\nabla f| \cdot |\nabla Ric|$, hence $\frac{\nabla f\cdot \nabla(\lambda_1+\lambda_2)}{f}$ is sufficiently small outside a compact set.

In all, we get
\ban
\Delta_f\frac{\lambda_1+\lambda_2}{f}\leq -0.9 \cdot\frac{\lambda_1+\lambda_2}{f}
\ean
outside a compact set $D$.

Suppose $\lambda_1+\lambda_2$ is not identically zero on $M\setminus D$,
then we can apply similar argument as \cite{Chow-Lu-Yang} or \cite{Munteanu-Wang4} to derive that
\ban
\frac{\lambda_1+\lambda_2}{f}\geq \frac{a}{f}
\ean
for some small positive $a$ outside a compact set. So $\lambda_1+\lambda_2\geq a$, which contradict with the fact that  $\lambda_1+\lambda_2$ tends to zero at infinity.

So $\lambda_1+\lambda_2\equiv 0$ on $M\setminus D$,  hence the function
\ban
G=tr(Ric^3)-\frac{1}{2}|Ric|^2,
\ean
is $0$ on $M\setminus D$.

\smallskip
Because $G$ is an analytic function, has to be zero, we obtain that $G\equiv 0$ on $M$. Moreover, the equation $0=\Delta_f R=R-2|Ric|^2$ implies that
\ban
G=&tr(Ric^3)-|Ric|^2+\frac{1}{4}R\\
=&\sum_{i=1}^4(\lambda_i-\frac{1}{2})^2 \lambda_i=0.
\ean
Finally we get $\lambda_1=\lambda_2\equiv 0$ and $\lambda_3=\lambda_4\equiv \frac{1}{2}$ due to $Ric\geq 0$. This implies the Ricci curvature has constant rank $2$. Therefore, any 4-dimensional shrinking gradient Ricci soliton with $R=1$ is  isometric to a finite quotient of $\mathbb{R}^2 \times \mathbb{S}^2$ by \cite{FR16}.
\qed

\end{document}